\documentclass[12pt, amssymb]{article}
\usepackage{amsmath}
\usepackage{amssymb}
\newenvironment{proof}{{\bf Proof}:}{\vskip 5mm }

\newtheorem{proposition}{Proposition}[subsection]
\newtheorem{lemma}[proposition]{Lemma}

\newtheorem{theorem}[proposition]{Theorem}

\newcommand{\llabel}[1]{\label{#1}}
\newcommand{\comment}[1]{}

\newcommand{\sr}{\rightarrow}
\newcommand{\zz}{{\bf Z}}
\newcommand{\hh}{{\bf{H}}}

\newcommand{\uu}{\underline}

\newcommand{\af}{{\bf A}^1}

\newcommand{\Ho}[3]{\underline{H}^{#1,#2}(#3)}
\newcommand{\motc}[3]{{H}^{#2,#3}(#1,\zz/2)}

\newcommand{\ep}{{\tau}}

\newcommand{\Hh}{{\mathbb{H}}}

\newcommand{\hiqa}{{\mathcal X}_{\uu{a}}}
\newcommand{\whiqa}{\widetilde{\mathcal X}_{\uu{a}}}
\newcommand{\mhiqa}{M({\mathcal X}_{\uu{a}})}
\newcommand{\ma}{{M}_{\uu{a}}}
\newcommand{\qa}{Q_{\uu{a}}}
\newcommand{\mqa}{M(Q_{\uu{a}})}

\newcommand{\K}{{\rm{K}}}

\newcommand{\Hom}{{\rm{H}\rm{o}\rm{m}}}

\newcommand{\la}{\langle}
\newcommand{\ra}{\rangle}
\newcommand{\lva}{\langle\!\langle}   
\newcommand{\rva}{\rangle\!\rangle}   

\newcommand{\grwk}{{\rm{G}}{\rm{r}}^*_{I^{\cdot}}(W(k))}

\newcommand{\grwkd}[1]{{\rm{G}}{\rm{r}}^{#1}_{I^{\cdot}}(W(k))}
\newcommand{\grwd}[2]{{\rm{G}}{\rm{r}}^{#2}_{I^{\cdot}}(W(#1))}

\newcommand{\nichego}[1]{}

\begin{document}
\begin{center}
{\Large\bf An exact sequence for $\K^M_{*}/2$ with applications to
quadratic forms}\\ D. Orlov\footnote{supported by the NSF grant DMS
97-29992 and the RFFI--99--01--01144} A. Vishik\footnote{supported by
the NSF grant DMS 97-29992} V. Voevodsky\footnote{Supported by the NSF
grants DMS-97-29992 and DMS-9901219 and The Ambrose Monell Foundation}
\end{center}
\vskip 4mm
\tableofcontents

\subsection{Introduction}
Let $k$ be a field of characteristics zero. For a sequence
$\uu{a}=(a_1,\dots,a_n)$ of invertible elements of $k$ consider the
homomorphism 
$$K_{*}^M(k)/2\sr K_{*+n}^M(k)/2$$
in Milnor's K-theory modulo elements divisible by $2$ defined by
the multiplication with the symbol corresponding to $\uu{a}$. The goal of
this paper is to construct a four-term exact sequence (\ref{mainseqf})
which provides information about the kernel and cokernel of this
homomorphism. 

The proof of our main theorem (Theorem \ref{mainf}) consists of two
independent parts. Let $Q_{\uu{a}}$ be the norm quadric defined by the
sequence $\uu{a}$ (see below). First, we use the techniques of
\cite{Vo4} to establish a four term exact sequence (\ref{mainseq})
relating the kernel and cokernel of the multiplication by $\uu{a}$
with Milnor's K-theory of the closed and the generic points of
$Q_{\uu{a}}$ respectively. This is done in the first section. Then,
using elementary geometric arguments, we show that the sequence can be
rewritten in its final form (\ref{mainseqf}) which involves only the
generic point and the closed points with residue fields of degree $2$.

As an application we establish, for fields of characteristics zero,
the validity of three conjectures in the theory of quadratic forms -
the Milnor conjecture on the structure of the Witt ring, the
Khan-Rost-Sujatha conjecture and the J-filtration conjecture. All
these results require only the first form of our exact sequence. Using
the final form of the sequence we also show that the kernel of
multiplication by $\uu{a}$ is generated, as a $K_*^M(k)$-module, by
its components of degree $\le 1$.

This paper is a natural extension of \cite{Vo4} and we feel free to
refer to the results of \cite{Vo4} without reproducing them here. Most
of the mathematics used in this paper was developed in the spring of
1996 when all three authors were at Harvard. In its present form the
paper was written while the authors were members of the Institute for
Advanced Study in Princeton. We would like to thank both institutions
for their support.

\subsection{An exact sequence for $\K^M_{*}/2$}
Let $\uu{a}=(a_1,..., a_n)$ be a sequence of elements of $k^*$.
Recall that the n--fold Pfister form $\lva a_1,..., a_n \rva$ is
defined as the tensor product
$$
\langle 1, -a_1\rangle\otimes \cdots\otimes \langle 1, -a_n\rangle
$$
where $\langle 1, -a_i\rangle$ is the norm form in the quadratic
extension $k(\sqrt{a_i})$. Denote by $\qa$ the projective quadric of
dimension $2^{n-1}-1$ defined by the form $q_{\uu{a}}=\lva a_1,...,
a_{n-1} \rva - \langle a_{n}\rangle$.  This quadric is called the
small Pfister quadric or the norm quadric associated with the symbol
$\uu{a}$. Denote by $k(Q_{\uu{a}})$ the function field of $Q_{\uu{a}}$
and by $(Q_{\uu{a}})_0$ the set of closed points of $Q_{\uu{a}}$. The
following result is the main theorem of the paper.
\begin{theorem}\label{main}
Let $k$ be a field of characteristic zero. Then for any sequence of
invertible elements $(a_1,\dots,a_n)$ the following sequence of
abelian groups is exact 
\begin{equation}\label{mainseq}
\mathop{\coprod}\limits_{x\in(\qa)_{(0)}}\K^M_{*}(k(x))/2
\stackrel{\rm{Tr_{k(x)/k}}}{\sr}
\K^M_{*}(k)/2\stackrel{\cdot\uu{a}}{\sr}
\K^{M}_{*+n}(k)/2\sr \K^{M}_{*+n}(k(\qa))/2
\end{equation}
\end{theorem}
The proof goes as follows. We first construct two exact sequences of
the form 
\begin{equation}\label{mainseq1}
0\sr K\sr \K^{M}_{*+n}(k)/2\sr \K^{M}_{*+n}(k(\qa))/2
\end{equation}
and
\begin{equation}\label{mainseq2}
\mathop{\coprod}\limits_{x\in(\qa)_{(0)}}\K^M_{*}(k(x))/2
\stackrel{\rm{Tr_{k(x)/k}}}{\sr}
\K^M_{*}(k)/2\sr I\sr 0
\end{equation}
and then construct an isomorphism $I\sr K$ such that the composition
$$\K^M_{*}(k)/2\sr I\sr K\sr \K^{M}_{*+n}(k)/2$$
is the multiplication by $\uu{a}$. 

Our construction of the sequence (\ref{mainseq1}) makes sense for any
smooth scheme $X$ and we shall do it in this generality. Recall that
we denote by $\check{C}(X)$ the simplicial scheme such that
$\check{C}(X)_n=X^{n+1}$ and faces and degeneracy morphisms are given
by partial projections and diagonal embeddings respectively. We will
use repeatedly the following lemma which is an immediate corollary of
\cite[Proposition 2.7]{Vo4} and \cite[Corollary 2.13]{Vo4}.
\begin{lemma}
\llabel{useful}
For any smooth scheme $X$ over $k$ and any $p\le q$ the homomorphism
$$H^{p,q}(Spec(k),\zz/2)\sr H^{p,q}(\check{C}(X),\zz/2)$$
defined by the canonical morphism $\check{C}(X)\sr Spec(k)$, is an
isomorphism.
\end{lemma}
\begin{proposition}
For any $n\ge 0$ there is an exact sequence of the form
\begin{equation}\label{s1}
0\sr H^{n,n-1}({\check{C}(X)},\zz/2)\sr
K_{n}^M(k)/2\sr K^M_{n}(k(X))/2
\end{equation}
\end{proposition}
\begin{proof}
The computation of motivic cohomology of weight $1$ shows that
$$\Hom(\zz/2, \zz/2(1))\cong H^{0,1}(Spec(k), \zz/2)\cong \zz/2$$
The nontrivial element $\ep: \zz/2\sr\zz/2(1)$ together with the
multiplication morphism $\zz(n-1)\otimes
\zz/2(1)\stackrel{\sim}{\sr}\zz/2(n)$ defines a morphism
$\ep:\zz/2(n-1)\sr \zz/2(n)$. The Beilinson-Lichtenbaum conjecture
implies immediately the following result.
\begin{lemma}
The morphism $\ep$ extends to a distinguished triangle in 
$DM^{eff}_{-}$ of the form
\begin{equation}\label{tens}
\zz/2(n-1)\stackrel{\cdot \ep}{\sr} \zz/2(n)\sr {\Ho{n}{n}{\zz/2}[-n]},
\end{equation}
where $\underline{H}^{n}(\zz/2(n))$ is the $n$\!--th cohomology sheaf of
the complex $\zz/2(n)$.
\end{lemma}
Consider the long sequence of morphisms in the triangulated category
of motives from the motive of $\check{C}(X)$ to the distinguished
triangle (\ref{tens}). It starts as
$$0\sr H^{n,n-1}(\check{C}(X), \zz/2)\sr H^{n,n}(\check{C}(X), \zz/2)\sr
H^{0}(\check{C}(X), \Ho{n}{n}{\zz/2})$$
By Lemma \ref{useful} there are isomorphisms
$$H^n(\check{C}(X),\zz/2(n))=H^{n,n}(Spec(k),\zz/2)=K^{M}_{n}(k)/2$$
On the other hand, since $\Ho{n}{n}{\zz/2}$ is a homotopy invariant
sheaf with transfers, we have an embedding
$$H^{0}(\check{C}(X), \Ho{n}{n}{\zz/2})\hookrightarrow
\Ho{n}{n}{\zz/2}(Spec(k(X)))$$
The right hand side is isomorphic to
$H^{n,n}(Spec(k(X)),\zz/2)=K_n^M(k(X))/2$. This completes the proof of
the proposition.
\end{proof}
Let us now construct the exact sequence (\ref{mainseq2}).  Denote the
standard simplicial scheme $\check{C}(\qa)$ by $\hiqa$.  Recall that
we have a distinguished triangle of the form
\begin{equation}\label{RMdec}
\mhiqa(2^{n-1}-1)[2^n-2]\stackrel{\varphi}{\to}\ma
\stackrel{\psi}{\to}\mhiqa\stackrel{\mu'}{\to}
\mhiqa(2^{n-1}-1)[2^n-1]
\end{equation}
where $M_{\uu{a}}$ is a direct summand of the motive of the quadric
$Q_{\uu{a}}$. Denote the composition
\begin{equation}\label{seq1}
\mhiqa\stackrel{\mu'}{\sr}\mhiqa(2^{n-1}-1)[2^n-1]
\stackrel{pr}{\sr}
\zz/2(2^{n-1}-1)[2^n-1]
\end{equation}
by $\mu\in\motc{\hiqa}{2^n-1}{2^{n-1}-1}$. By Lemma \ref{useful} we
have
$$H^{i,i}({\cal X}_{\uu{a}},\zz/2)=H^{i,i}(Spec(k),\zz/2)=K_i^M(k)/2$$
Therefore, multiplication with $\mu$ defines a homomorphism
$$\K^M_{i}(k)/2\stackrel{\cdot\mu}{\sr}
\motc{\hiqa}{i+2^n-1}{i+2^{n-1}-1}.$$
\begin{proposition}\llabel{seq2}
The sequence
\begin{equation}\label{s2}
\mathop{\coprod}\limits_{x\in(\qa)_{(0)}}\K^M_{i}(k(x))/2
\stackrel{\rm{Tr_{k(x)/k}}}{\sr}
\K^M_{i}(k)/2\stackrel{\cdot\mu}{\sr}
\motc{\hiqa}{i+2^n-1}{i+2^{n-1}-1}\sr 0
\end{equation}
is exact.
\end{proposition}
\begin{proof}
Let us consider morphisms in the triangulated category of motives from
the distinguished triangle (\ref{RMdec}) to the object
$\zz/2(i+2^{n-1}-1)[i+2^n-1]$.  By definition, $\ma$ is a direct
summand of the motive of the smooth projective variety $\qa$ of
dimension $2^{n-1}-1$, therefore, the group
$\motc{\ma}{i+2^{n}-1}{i+2^{n-1}-1}$ is trivial by \cite[Corollary
2.3]{Vo4}.  Using this fact, we obtain the following exact sequence:
\begin{equation}
\llabel{eq9.28.1}
\motc{\ma}{i+2^{n}-2}{i+2^{n-1}-1}\stackrel{\varphi^*}{\sr}
\motc{\hiqa}{i}{i}\stackrel{{\mu'}^*}{\sr}
\end{equation}
$$
\sr\motc{\hiqa}{i+2^{n}-1}{i+2^{n-1}-1}{\sr}0
$$
By definition (see \cite[p.43]{Vo4}) the morphism $\varphi$ is given
by the composition
\begin{equation}
\llabel{eq9.28.2}
\mhiqa(2^{n-1}-1)[2^n-2]\stackrel{pr}{\sr}
\zz(2^{n-1}-1)[2^n-2]\sr M_{\uu{a}}
\end{equation} 
and the composition of the second arrow with the canonical embedding
$M_{\uu{a}}\sr M(Q_{\uu{a}})$ is the fundamental cycle map 
$$\zz(2^{n-1}-1)[2^n-2]\sr M(Q_{\uu{a}})$$
which corresponds to the fundamental cycle on $\qa$ under the
isomorphism
$$\Hom(\zz(2^{n-1}-1)[2^n-2],\mqa)={\rm{CH}}_{2^{n-1}-1}(\qa)\cong \zz$$
(see \cite[Theorem 4.5]{Vo4}). On the other hand by Lemma \ref{useful}
the homomorphism 
$$\motc{Spec(k)}{i}{i}\sr \motc{\hiqa}{i}{i}$$
defined by the first arrow in (\ref{eq9.28.2}) is an isomorphism.
This implies immediately that the exact sequence (\ref{eq9.28.1})
defines an exact sequence of the form
\begin{equation}
\llabel{eq9.28.3}
\motc{Q_{\uu{a}}}{i+2^{n}-2}{i+2^{n-1}-1}\stackrel{\varphi^*}{\sr}
\motc{Spec(k)}{i}{i}\stackrel{{\mu'}^*}{\sr}
\end{equation}
$$
\sr\motc{\hiqa}{i+2^{n}-1}{i+2^{n-1}-1}{\sr}0
$$
By \cite[Corollary 2.4]{Vo4} there is an isomorphism
$$
\motc{\qa}{i+2^{n}-2}{i+2^{n-1}-1}\cong 
H^{2^{n-1}-1}(\qa, \underline{K}^{M}_{i+ 2^{n-1} -1}/2)
$$ 
The Gersten resolution for the sheaf $\uu{K}_m^M/2$ (see, for example,
\cite{R1}) shows that the group $H^{2^{n-1}-1}(\qa,
\underline{K}^{M}_{i+ 2^{n-1} -1}/2)$ can be identified with the
cokernel of the map:
$$
\mathop{\coprod}\limits_{y\in(\qa)_{(1)}}{\K}^M_{i+1}(k(y))/2
\stackrel{\partial}{\sr}
\mathop{\coprod}\limits_{x\in(\qa)_{(0)}}{\K}^M_{i}(k(x))/2,
$$
and the map $\motc{\qa}{i+2^{n}-2}{i+2^{n-1}-1}{\to}
\motc{{Spec(k)}}{i}{i}$ defined by the fundamental cycle corresponds
in this description to the map
$$
\mathop{\coprod}\limits_{x\in(\qa)_{(0)}}{\K}^M_{i}(k(x))/2
\stackrel{\rm{Tr_{k(x)/k}}}{\sr}
{\K}^M_{i}(k)/2=\motc{{Spec}(k)}{i}{i}
$$ 
This finishes the proof of Proposition \ref{seq2}.
\end{proof}
We are going to show now that the map
$\K^{M}_{*}(k)/2\stackrel{\uu{\alpha}}{\sr}\K^{M}_{*+n}(k)/2$ glues
the exact sequences (\ref{s1}) and (\ref{s2}) in one. Denote by
$\Hh^{i}(\hiqa)$ the direct sum $\oplus_{m} \motc{\hiqa}{m+i}{m}$. It
has a natural structure of a graded module over the ring
$\K^{M}_{*}(k)/2$ and one can easily see that the sequences (\ref{s1})
and (\ref{s2}) define sequences of
$\K^{M}_{*}(k)/2$\!--modules of the form
\begin{eqnarray}
&0\sr\Hh^{1}(\hiqa)\sr\K^{M}_{*}(k)/2\sr \K^{M}_{*}(k(\qa))/2\label{first}&\\
&\mathop{\coprod}\limits_{x\in(\qa)_{(0)}}\K^M_{*}(k(x))/2
\stackrel{\rm{Tr_{k(x)/k}}}{\sr}
\K^M_{*}(k)/2\stackrel{\cdot\mu}{\sr}
\Hh^{2^{n-1}}(\hiqa)\sr 0 &\label{second}
\end{eqnarray}
Consider cohomological operations
$$
Q_{i}: \motc{-}{\bullet}{*}\sr \motc{-}{\bullet+2^{i+1}-1}{*+2^i -1}
$$
introduced in \cite[p.32]{Vo4}. The composition $Q_{n-2}\cdots Q_0$
defines a homomorphism of graded abelian groups $d:\Hh^{1}(\hiqa)\sr
\Hh^{2^{n-1}}(\hiqa)$ and \cite[Theorem 3.17(2)]{Vo4} together with
the fact that $H^{p,q}(Spec(k),\zz/2)=0$ for $p>q$ implies that $d$ is
a homomorphism of $K_*^M(k)/2$-modules. We are going to show that $d$
is an isomorphism and that the composition
\begin{equation}
\llabel{compeq}
\K^M_{*}(k)/2\stackrel{\cdot\mu}{\sr}
\Hh^{2^{n-1}}(\hiqa)\stackrel{d^{-1}}{\sr}\Hh^{1}(\hiqa)\sr\K^{M}_{*}(k)/2
\end{equation}
is the multiplication with $\uu{a}$.
\begin{lemma} 
\llabel{inj}
The homomorphism $d$ is injective.
\end{lemma}
\begin{proof}
We have to show that the composition of operations 
$$Q_{n-2}\dots Q_0:\motc{\hiqa}{*+n}{*+n-1}{\sr}
\motc{\hiqa}{*+2{n}-1}{*+2^{n-1} -1}$$ 
is injective. Let $\whiqa$ be the simplicial cone of the morphism
$\hiqa\sr Spec(k)$ which we consider as a pointed simplicial
scheme. The long exact sequence of cohomology defined by the
cofibration sequence 
\begin{equation}
\label{cofseq}
(\hiqa)_+\sr Spec(k)_+\sr \whiqa\sr \Sigma^1_s((\hiqa)_+)
\end{equation}
together with the fact that $H^{p,q}(Spec(k),\zz/2)=0$ for $p>q$ shows
that for $p>q+1$ we have a natural isomorphism
$H^{p,q}(\whiqa,\zz/2)=H^{p-1,q}(\hiqa,\zz/2)$ compatible with the
actions of cohomological operations. Therefore, it is sufficient to
prove injectivity of the composition $Q_{n-2}\dots Q_0$ on motivic
cohomology groups of the form $\motc{\whiqa}{*+n+1}{*+n-1}$. To show
that $Q_{n-2}\cdots Q_0$ is a monomorphism it is sufficient to check
that the operation $Q_i$ acts monomorphically on the group
$$\motc{\whiqa}{*+n-i+2^{i+1}-1}{*+n-i-2^{i}-2}$$
for all $i=0,..., n-2$. For any $i\le n-1$ we have $ker(Q_i)=Im(Q_i)$
by \cite[Theorem 3.25]{Vo4} and \cite[Lemma 4.11]{Vo4}. Therefore,
the kernel of $Q_i$ on our group is the image of
$\motc{\whiqa}{*+n-i)}{*+n-i-1}$. On the other hand, the cofibration
sequence (\ref{cofseq}) together with Lemma \ref{useful} implies that
for $p\le q+1$ we have $H^{p,q}(\whiqa,\zz/2)=0$ which proves the
lemma.
\end{proof}
Denote by $\gamma$ the element of $\motc{\hiqa}{n}{n-1}$ which
corresponds to the symbol $\uu{a}$ under the embedding into
$\K^{M}_{n}(k)/2$ (sequence (\ref{s1})). To prove that $d$ is
surjective and that the composition (\ref{compeq}) is multiplication with
$\uu{a}$ we use the following lemma.
\begin{lemma}
\llabel{double}
The composition
$\K^{M}_{*}(k)/2\stackrel{\cdot \gamma}{\sr} \Hh^1(\hiqa)
\stackrel{d}{\sr}\Hh^{2^{n-1}}(\hiqa)$
coincides with the multiplication by $\mu$.
\end{lemma}
\begin{proof}
Since our maps are homomorphisms of $K_*^M(k)$-modules it is
sufficient to verify that the cohomological operation $d$ sends
$\gamma\in \motc{\hiqa}{n}{n-1}$ to $\mu\in
\motc{\hiqa}{2^n-1}{2^{n-1}-1}$. By Lemma \ref{inj} $d$ is
injective. Therefore, the element $d(\gamma)$ iz nonzero.  On the
other hand, sequence (\ref{s2}) shows that 
$$
\motc{\hiqa}{2^n-1}{2^{n-1}-1}\cong \K^{M}_{0}(k)/2\cong \zz/2 
$$
and $\mu$ is a generator of this group. Therefore, $d(\gamma)=\mu$.
\end{proof}
\begin{lemma}
\llabel{sur}
The homomorphism $d$ is surjective.
\end{lemma}
\begin{proof}
Follows immediately from Lemma \ref{double} and surjectivity of
multiplication by $\mu$ (Proposition \ref{seq2}).
\end{proof}
\begin{lemma}
\llabel{comp}
The composition (\ref{compeq}) is the multiplication with $\uu{a}$.
\end{lemma}
\begin{proof}
Since all the maps in (\ref{compeq}) are morphisms of $K^M_*(k)$-modules,
it is sufficient to check the condition for the generator $1\in\K^M_0(k)/2$.
And the later follows from Lemma \ref{double} and the definition of $\gamma$.
\end{proof}
This finishes the proof of Theorem \ref{main}.
\vskip 3mm

%
\noindent
The following statement, which is easily deduced  from 
the exact sequence (\ref{mainseq}), is the key to many
applications.

Let $E/k$ be a field. For any element $h\in \K^M_n(k)$ 
denote by $h|_{E}$, as usually, the restriction of $h$ on $E$, i.e.
image of $h$ under the natural morphism $\K^M_n(k)\to \K^M_n(E)$.
\begin{theorem}\label{conseq}
For any field $k$ and any nonzero $h\in\K^M_n(k)/2$ there exist a field 
$E/k$ and a pure symbol $\alpha=\{a_1,\dots,a_n\}\in\K^M_n(k)/2$ such that
$h|_{E}=\alpha|_{E}$ is a {\bf nonzero} pure symbol of $\K^M_n(E)/2$.
\end{theorem}
\begin{proof}
Let $h=\alpha_1+\dots+\alpha_l$, where $\alpha_i$ are pure symbols
corresponding to sequences $\uu{a}_i=(a_{1i},\dots,a_{ni})$.  Let
$Q_{\uu{a}_i}$ be the norm quadric corresponding to the
symbol $\alpha_i$.  For any $0<i\le l$ denote by $E_i$ the field
$k(Q_{\uu{a}_1}\times\dots\times Q_{\uu{a}_i})$.  It is clear that
$h|_{E_l}=0$. Let us fix $i$ such that $h|_{E_{i+1}}=0$ and $h|_{E_i}$
is a nonzero element. Then $h|_{E_i}$ belongs to
$${ker}(\K^M_n(E_i)/2\to\K^M_n(E_{i+1})/2)$$
By Theorem \ref{main}, the kernel is covered by $K^M_0(E_{i})\cong
\zz/2$ and is generated by $\alpha_{i+1}|_{E_{i}}$.  Thus, we have
$\alpha_{i+1}|_{E_{i}}=h|_{E_i}\ne 0$.
\end{proof}

\subsection{Reduction to points of degree $2$}
In this section we prove the following result.
\begin{theorem}
\llabel{reduction} Let $k$ be a field such that $char(k)\ne 2$ and $Q$
be a smooth quadric over $k$. Let $Q_{(0)}$ be the set of closed
points of $Q$ and $Q_{(0,\le 2)}$ the subset in $Q_{(0)}$ of points
$x$ such that $[k_x:k]\le 2$. Then, for any $n\ge 0$, the image of the
map
\begin{equation}
\llabel{eq15}
\oplus tr_{k_x/k}:\oplus_{x\in Q_{(0)}}K_n^M(k_x)\sr K_n^M(k)
\end{equation}
coincides with the image of the map
\begin{equation}
\llabel{smalltr}
\oplus tr_{k_x/k}:\oplus_{x\in Q_{(0,\le 2)}}K_n^M(k_x)\sr
K_n^M(k)
\end{equation}
\end{theorem}
Combining Theorem \ref{main} with Theorem \ref{reduction} we get the
following result.
\begin{theorem}
\llabel{mainf}
\llabel{cor1}
Let $k$ be a field of characteristic zero and $\uu{a}=(a_1,\dots, a_n)$
a sequence of invertible elements of $k$.  Then the sequence
\begin{equation}
\llabel{mainseqf} 
\oplus_{x\in (Q_{\uu{a}})_{(0,\le 2)}}K_i^M(k_x)/2\sr
K_{i}^M(k)/2\stackrel{\uu{a}}{\sr}K_{i+n}^M(k)/2\sr
K^M_{i+n}(k(Q_{\uu{a}}))/2
\end{equation}
is exact.
\end{theorem}
Theorem \ref{cor1} together with the well known result of Bass and
Tate (see \cite[Corollary 5.3]{BT}) implies the following.
\begin{theorem}
\llabel{cor2} Let $k$ be a field of characteristic zero and
$\uu{a}=(a_1,\dots, a_n)$ a sequence of invertible elements of $k$
such that the corresponding elements of $K_{n}^M(k)/2$ is not zero.
Then the kernel of the homomorphism
$K_{*}^M(k)/2\stackrel{\uu{a}}{\sr}K_{*+n}^M(k)/2$ is generated, as a
module over $K_{*}^M(k)$, by the kernel of the homomorphism
$K_1^M(k)/2\sr K_{1+n}^M(k)/2$.
\end{theorem}
Let us start the proof of Theorem \ref{reduction} with the following
two lemmas.
\begin{lemma}
\llabel{simple} Let $E$ be an extension of $k$ of degree $n$ and $V$ a
$k$-linear subspace in $E$ such that $2dim(V)>n$. Then, for any $n>0$,
$K_n^M(E)$ is generated, as an abelian group, by elements of the
form $(x_1,\dots, x_n)$ where all $x_i$'s are in $V$.
\end{lemma}
\begin{proof}
It is sufficient to prove the statement for $n=1$. Let $x$ be an
invertible element of $E$. Since $2dim(V)>dim_kE$ we have $V\cap xV\ne
0$. Therefore $x$ is a quotient of two elements of $V\cap E^*$.
\end{proof}
\begin{lemma}
\llabel{simple2} Let $k$ be an infinite field and $p$ a closed separable point
in ${\bf P}^n_k$, $n\ge 2$ of degree $m$. Then there exists a rational
curve $C$ of degree $m-1$ such that $p\in C$ and $C$ is either nonsingular,
or has one rational singular point.
\end{lemma}
\begin{proof}
We may assume that $p$ lies in ${\bf A}^n\subset {\bf P}^n$. Then
there exists a linear function $x_1$ on ${\bf A}^n$ such that the map
of the residue fields $k_{x_1(p)}\sr k_p$ is an isomorphism. Let
$(x_1,\dots,x_n)$ be a coordinate system starting with $x_1$. Since
the restriction of $x_1$ to $p$ is an isomorphism the inverse gives a
collection of regular functions $\bar{x}_2,\dots,\bar{x}_n$ on
$x_1(p)\subset \af$. Each of this functions has a representative $f_i$
in $k[x_1]$ of degree at most $m-1$. The projective closure of the
affine curve given by the equations $x_i=f_i(x_1)$, $i=2,\dots,n$
satisfies the conditions of the lemma.
\end{proof}
Let $Q$ be any quadric over $k$. If $Q$ has a rational point (or even
a point of odd degree, which is the same by Springer's theorem,
\cite[VII,Theorem 2.3]{Lam}), then Theorem \ref{reduction} for $Q$
holds for obvious reasons. Therefore we may assume that $Q$ has no
points of odd degree. It is well known (see e.g. \cite[Th.2.3.8,
p.39]{Sch}) that any smooth quadric of dimension $>0$ over a finite
field of odd characteristic has a rational point. Since the statement
of the theorem is obvious for $dim(Q)=0$ we may assume that $k$ is
infinite. By the theorem of Springer, for finite extension of odd degree 
$E/F$, the quadric $Q_F$ is isotropic if and only if $Q_E$ is. Hence, we
can assume that $E/k$ is separable.

Let $e$ be a point on $Q$ with the residue field $E$. We have to show
that the image of the transfer map $K_n^M(E)\sr K_n^M(k)$ lies in the
image of the map (\ref{smalltr}). We proceed by induction on $d$ where
$2d=[E:k]$. If $d=1$ there is nothing to prove. Assume by induction
that for any closed point $f$ of $Q$ such that $[k_f:k]<2d$ the image
of the transfer map $K_n^M(k_f)\sr K_n^M(k)$ lies in the image of
(\ref{smalltr}).  

If $dim(Q)=0$ our statement is obvious. Consider the case of a conic
$dim(Q)=1$. Let $D$ be any effective divisor on $Q$ of degree
$2d-2$. Denote by $h^0(D)$ the linear space $H^0(Q,{\cal O}(D))$ which
can be identified with the space of rational functions $f$ such that
$D+(f)$ is effective.  Evaluating elements of $h^0(D)$ on $e$ we get a
homomorphism $h^0(D)\sr E$ which is injective since $deg(D)<2d$. By
the Riemann-Roch theorem we have $dim(h^0(D))=2d-1$ and therefore, by Lemma
\ref{simple}, $K_n^M(E)$ is generated by elements of the form
$\{f_1(e),\dots,f_n(e)\}$ where $f_i\in h^0(D)$. Let now $D'$ be an
effective divisor on $Q$ of degree $2$ (it exists since $Q$ is a
conic). Using again the Riemann-Roch theorem we see that $dim(|e-D'|)>0$
i.e. that there exists a rational function $f$ with a simple pole in
$e$ and a zero in $D'$. In particular, the degrees of all the points
where $f$ has singularities other than $e$ is strictly less than
$2d$. Consider the symbol $\{f_1,\dots,f_n,f\}\in
K_{n+1}^M(k(Q))$. Let
$$\partial:K_{n+1}^M(k(Q))\sr \oplus_{x\in Q_{(0)}} K_n^M(k_x)$$
be the residue homomorphism. By \cite[]{} its composition with
(\ref{eq15}) is zero. On the other hand we have 
$$\partial(\{f_1,\dots,f_n,f\})=\{f_1(e),\dots,f_n(e)\}+u$$
where $u$ is a sum of symbols concentrated in the singular points of
$f_1,\dots, f_n$ and singular points of $f$ other than $e$. Therefore,
by our construction $u$ belongs to $\oplus_{x\in
Q_{(0),<2d}}K_n^M(k_x)$ and we conclude that
$tr_{E/k}\{f_1(e),\dots,f_n(e)\}$ lies in the image of (\ref{smalltr})
by induction.

Let now $Q$ be a quadric in ${\bf P}^n$ where $n\ge 3$. Let $c$ be a
rational point of ${\bf P}^n$ outside $Q$ and $\pi:Q\sr {\bf P}^{n-1}$
be the projection with the center in $c$. The ramification locus of
$\pi$ is a quadric on ${\bf P}^{n-1}$ which has no rational points.

Assume first that there exists $c$ such that the degree of $\pi(e)$ is
$d$. Then, by Lemma \ref{simple2}, we can find a (singular) rational
curve $C'$ in ${\bf P}^{n-1}$ of degree $d-1$ which contains
$\pi(e)$. Consider the curve $C=\pi^{-1}(C')\subset Q$. Let
$\tilde{C},\tilde{C'}$ be the normalizations of $C$ and $C'$ and
$\tilde{\pi}:\tilde{C}\sr \tilde{C'}$ the morphism corresponding to
$\pi$. Since $deg(e)=2d$ and $deg(\pi(e))=d$ the point $e$ does not
belong to the ramification locus of $\pi:Q\sr {\bf P}^{n-1}$. This
implies that $e$ lifts to a point $\tilde{e}$ of $\tilde{C}$ of degree
$2d$ and that
$\tilde{e}=\tilde{\pi}^{-1}(\tilde{\pi}(\tilde{e}))$. Since the
ramification locus of $\pi$ has no rational points the singular point
of $C'$ is unramified. This implies that $\tilde{\pi}$ is ramified in
$\le 2(d-1)$ points and, therefore, $\tilde{C}$ is a hyperelliptic
curve of genus less or equal to $d-2$. Let $D$ be an effective divisor
on $\tilde{C}$ of degree $2d-2$. By the Riemann-Roch theorem we have
$dim(h^0(D))\ge d+1$. On the other hand, since $deg(D)<2d$, the
homomorphism $h^0(D)\sr E$ defined by evaluation at $\tilde{e}$ is
injective. Therefore, by Lemma \ref{simple} any element in $K_n^M(E)$
is of the form $\{f_1(\tilde{e}),\dots,f_n(\tilde{e})\}$ for $f_i\in
h^0(D)$.  Let $D'$ be an effective divisor on $\tilde{C}$ of degree
$2$. By the Riemann-Roch theorem we have $dim(h^0(\tilde{e}-D'))\ge
d+1>0$. Therefore, there exists a rational function $f$ with simple
pole in $\tilde{e}$ and such that all its other singularities are
located in points of degree $<2d$. We can conclude now that
$tr_{E/k}\{f_1(\tilde{e}),\dots,f_n(\tilde{e})\}$ belongs to the image
of (\ref{smalltr}) in the same way as we did in the case of
$dim(Q)=1$.

Consider now the general case - we still assume that $n\ge 3$ but not
that we can find a center of projection $c$ such that
$deg(\pi(e))=d$. Taking a general $c$ we may assume that
$deg(\pi(e))=2d$ and that $e$ does not belong to the ramification
locus of $\pi$. By Lemma \ref{simple2} we can find a rational curve
$C'$ in ${\bf P}^{n-1}$ of degree $2d-1$ which contains
$\pi(e)$. Consider the curve $C=\pi^{-1}(C')\subset Q$. Let
$\tilde{C},\tilde{C'}$ be the normalizations of $C$ and $C'$ and
$\tilde{\pi}:\tilde{C}\sr \tilde{C'}$ the morphism corresponding to
$\pi$. Since the point $e$ does not belong to the ramification locus
of $\pi$ it lifts to a point $\tilde{e}$ of $\tilde{C}$ of degree
$2d$. Since the ramification locus of $\pi$ does not have rational
points and the only singular point of $C'$ is rational, $\tilde{\pi}$
is ramified in no more than $2(2d-1)$ points and, therefore,
$\tilde{C}$ is a hyperelliptic curve of genus $\le 2d-2$. Let $D$ be
an effective divisor on $\tilde{C}'={\bf P}^1$ of degree $d$.  We have
$dim(h^0(D))=d+1$ and since $\tilde{\pi}(\tilde{e})$ has degree $2d$
the evaluation at $\tilde{\pi}(\tilde{e})$ gives an injective
homomorphism $h^0(D)\sr E$. By Lemma \ref{simple}, we conclude that
any element of $K_n^M(E)$ is a linear combination of elements of the
form $\{f_1(\tilde{\pi}(\tilde{e})),\dots,
f_n(\tilde{\pi}(\tilde{e}))\}$ where $f_i$ are in $h^0(D)$.  Let $D'$
be an effective divisor of degree $2$ on $\tilde{C}$. By the
Riemann-Roch theorem we have $dim(h^0(\tilde{e}-D'))\ge 1$ i.e. there exists
a rational function $f$ with a simple pole in $\tilde{e}$ and a zero
in $D'$. If $d>1$ then all the singular points of $f$, but $\tilde{e}$,
are of degree $<2d$ and by the same reasoning as in the previous two
cases we conclude that $tr_{E/k}(\{f_1(\tilde{\pi}(\tilde{e})),\dots,
f_n(\tilde{\pi}(\tilde{e}))\})$ is a linear combination of the form
$$\sum_{x\in \tilde{C}_{(0),<2d}} tr_{k_{x}/k}(u_{x})+\sum_{i,y\in
(f_i\circ\tilde{\pi})} tr_{k_{y}/k}(v_{i,y})$$
Summands of the first type are in the image of (\ref{smalltr}) by the
inductive assumption. The fact that summands of the second type are in
the image of (\ref{smalltr}) follows from the case $deg(\pi(e))=d$
considered above.

\subsection{Some applications}
\subsubsection{Milnor's Conjecture on quadratic forms.}
As the first corollary of Theorem \ref{conseq} 
we get {\it Milnor's Conjecture on
quadratic forms}.
 
As usually, we denote by $W(k)$ the {\it Witt ring} of quadratic forms
over $k$, and by $I\subset W(k)$ the ideal of even-dimensional
forms. The filtration $W(k)\supset I\supset I^2\supset\dots\supset
I^n\supset\dots$ by the powers of $I$ is called the $I$-filtration on
$W$. We denote the associated graded ring by $\grwk$. Consider the map 
$$\K^M_1(k)/2=k^*/(k^*)^2\stackrel{\varphi_1}{\sr}\grwkd{1}$$
which sends $\{a\}$ to $\la 1,-a\ra$. Since $(\la 1,-a\ra+\la
1,-b\ra-\la 1,-ab\ra)\in I^2$ it is a group-homomorphism and one can
easily see that it is an isomorphism. For any $a\in k^*\backslash 1$,
the form $\lva a,1-a\rva$ is hyperbolic and, therefore, the
isomorphism $\varphi_1$ can be extended to a ring homomorphism
$\varphi:\K^M_*(k)/2\to\grwk$. Since $\grwk$ is generated by the
first-degree component $\varphi$ is surjective.  The {\it Milnor
Conjecture on quadratic forms} states that $\varphi$ is an isomorphism
i.e. that it is injective.  It was proven in degree 2 by J.Milnor
\cite{Mil}, in degree 3 by M.Rost \cite{R0} and A.Merkurjev-A.Suslin
\cite{MS}, and in degree 4 by M.Rost.  Moreover, R.Elman and T.Y.Lam
\cite{EL} proved that the map $\varphi$ is injective on {\it pure
symbols}.
\begin{theorem}\label{Mcon} Let $k$ be a field of characteristic
zero. Then, the natural map $\varphi:\K^M_*(k)/2\to\grwk$ is an
isomorphism.
\end{theorem}
\begin{proof}
We already know that $\varphi$ is surjective. Let $h\ne 0$ be an
element of $\K^M_n(k)/2$. By Theorem \ref{conseq} there exists
a field extension $E/k$ such that $h|_E$ is a nonzero pure symbol. By
a result of R.Elman and T.Y.Lam (\cite{EL}) the map $\varphi$
is injective on {\it pure symbols}.  Hence $\varphi(h|_{E})$ is
a nonzero element of $\grwd{E}{n}$.  Since the morphism $\varphi$ is
compatible with field extensions, the element $\varphi(h)\in
\grwd{k}{n}$ is also nonzero.  Therefore, $\varphi$ is injective.
\end{proof}

\subsubsection{Kahn-Rost-Sujatha Conjecture}
In \cite{KRS} B.Kahn, M.Rost and R.Sujatha proved that for any quadric
$Q$ of dimension $m$ the $ker(\K^M_i(k)/2\to\K^M_i(k(Q)))$ is trivial
for any $i<log_2(m+2)$, if $i\leq 4$ (actually, in \cite{KRS} authors
worked with $H^i_{et}(k, \zz/2)$ instead of $\K^M_i(k)/2$, but because
of \cite{Vo4} we can use $\K^M_i(k)/2$ here).  The authors also
conjectured
\footnote{only in the original version of the paper} (among other
things) that the same is true without the restriction $i\leq 4$. The
following result proves this conjecture.
\begin{theorem}\label{KRScon}
Let $Q$ be an $m$-dimensional quadric over a field $k$ of
characteristic zero. Then $ker(\K^M_i(k)/2\to\K^M_i(k(Q))/2)$ is
trivial for any $i<log_2(m+2)$.
\end{theorem}
\begin{proof}
Denote by $q$ a quadratic form which defines the quadric $Q$.  Assume
that $h$ is a nonzero element of $ker(\K^M_i(k)/2\to\K^M_i(k(Q))/2)$.
Using Theorem \ref{conseq} we can find an extension $E/k$ such that
$h|_E$ is a nonzero pure symbol of the form $\uu{a}=\{ a_1 ,...,
a_n\}$.  Then, since $h|_{E(Q)}=0$, the corresponding Pfister quadric
$\qa/E$ becomes hyperbolic over $E(Q)$.  Since $\qa|_{E(Q)}$ is
hyperbolic the form $t\cdot q|_E$ is isomorphic to a subform of the
Pfister form $\lva a_1 ,..., a_n \rva$ for some coefficient $t\in
E^*$ by \cite[ch.4, Theorem 5.4]{Sch}. In particular,
$m+2=\dim(Q)+2=dim(q)\le 2^i$. Therefore, we have $i\ge log_2(m+2)$.
\end{proof}

\subsubsection{$J$\!-- filtration conjecture.}

Together with the $I$\!-- filtration on $W(k)$ we can consider the
following so-called $J$\!-- filtration.  Let $x\in W(k)$ be an
element, $q$ an anisotropic quadratic form which represents $x$ and
$Q$ the corresponding projective quadric.  Since $Q$ has a point
over the field $k(Q)$, we have a decomposition of the form
$$
q|_{k(Q)}=q_1\perp \underset{i_1(q)}%
{\underbrace{\hh\perp...\perp\hh}},
$$ 
where $q_1$ is an anisotropic form over $k(Q)$, and $\hh$ is the
elementary hyperbolic form.  The number $i_1(q)$ is called the first
{\it higher Witt index} of $q$.  In the same way we can decompose
$q_1|_{k(Q)(Q_1)}$ etc. obtaining a sequence of quadratic forms $q,
q_1, ..., q_{s-1}$, where each $q_i$ is an anisotropic form defined
over $k(Q)...(Q_{i-1})$, and
$$
q_{s-1}|_{k(Q)...(Q_{s-1})}=\underset{i_{s}(q)}%
{\underbrace{\hh\perp...\perp\hh}}
$$ 
is a hyperbolic form.  By \cite[Theorem 5.8]{Kn} (see also \cite[ch
4.,Theorem 5.4]{Sch}), any quadratic form $q'$ over a field $E$, such
that $q' |_{E(Q')}$ is hyperbolic, is proportional to some Pfister
form.  This implies that the form $q_{s-1}$ is proportional to an
n--fold {\it Pfister form} $\lva a_1,\dots,a_n\rva$, where
$\{a_1,...,a_n\}\in \K^M_{n}(k(Q)...(Q_{s-2}))/2$.  This procedure
defines for any element $x\in W(k)$ a natural number $n$ which we will
call the degree of $x$.

Let us define $J_n(W(k))$ as the subset of $W(k)$ consisting of all
elements of {\it degree} $\ge n$.  It can be easily checked that
$I^n\subseteq J_n$.  It was conjectured in \cite[Question 6.7]{Kn} and
in \cite{Sch} that the $J_{\cdot}$ coincides with the $I^{\cdot}$.
The following theorem proves this conjecture.
\begin{theorem}
$J_{n}=I^{n}$. 
\end{theorem}
\begin{proof}
Let $x$ be an element of $J_n(W(k))$ which is represented by a
quadratic form $q$.  As above we have a sequence of quadrics $Q, Q_1,
..., Q_{s-1}$ such that $q|_{k(Q)(Q_1)...(Q_{s-1})}$ is
hyperbolic. This means that $x$ goes to $0$ under the natural map from
$W(k)$ to $W(k(Q)(Q_1)...(Q_{s-1}))$.

All quadrics $Q,Q_1,\dots,Q_{s-1}$ have dimensions $\geq
2^n-2>2^{n-1}-2$.  By Theorem \ref{KRScon}, for any $0\leq i\leq n-1$,
the kernel ${ker}(\K^M_i(k)/2\to\K^M_i(k(Q)...(Q_{s-1})))$ is trivial.
Therefore, applying the Milnor conjecture (Theorem \ref{Mcon}), we
conclude that the map
$$
\grwkd{i}\sr {{\rm{G}}{\rm{r}}^{i}_{I^{\cdot}}(W(k(Q)...(Q_{s-1})))} 
$$
is a monomorphism for all $0\le i\le n-1$. Therefore the map
$$
W(k)/I^n(W(k))\sr W(k(Q)...(Q_{s-1}))/I^n(W(k(Q)...(Q_{s-1})))
$$ 
is a monomorphism as well. Therefore, $x$ belongs to $I^n(W(k))$.
\end{proof}


\begin{thebibliography}{11}
%
\bibitem{BT}
H. Bass and J. Tate, \ {\it The {M}ilnor ring of a global field},
Lecture Notes in Math. {\bf 342} (1973), 349-446 

\bibitem{EL}
   Richard Elman, T.Y.Lam, \ {\it Pfister forms and K-theory of fields},
   J.Algebra {\bf 23} (1972), 181-213.
%
\bibitem{Lam}
   T.Y.Lam,  \ {\it Algebraic theory of quadratic forms},
   Benjamin, Readings, Mass., 1973.
%
\bibitem{Kn}
   M.Knebusch, \ 
   {\it Generic splitting of quadratic forms, I.},
   Proc. London Math. Soc. {\bf 33} (1976), 67-93.
%
\bibitem{KRS}
   B.Kahn, M.Rost, and R.J.Sujatha, \ 
   {\it Unramified cohomology of quadrics, I.},
   Amer.J.Math. {\bf 120} (1998), no.4, 841-891.
%
\bibitem{Mil}
   J.Milnor, \ {\it Algebraic K - theory and quadratic forms},
    Invent. Math. {\bf 9} (1969/70), 318-344.
%
\bibitem{MS}
   A.S.Merkurjev, A.A.Suslin, \ 
    {\it The norm-residue homomorphism of degree $3$} (in Russian),
    Izv. Akad. Nauk SSSR {\bf 54} (1990), 339-356.
    English translation: Math. USSR Izv. {\bf 36} (1991), 349-368.
%
\bibitem{R0}
   M.Rost, \ {\it Hilbert theorem $90$ for $K^M_3$ for degree-two extensions},
   Preprint, Regensburg, 1986.
%
\bibitem{R}
   M.Rost, \ {\it Some new results on the Chow-groups of quadrics},
    Preprint, Regensburg, 1990.
%
\bibitem{R1}
  M.Rost, \ {\it Chow groups with coefficients}, Doc. Math. 1 (1996), 319-393.
   MR 98a:14006, Zbl 864.14002.
%
\bibitem{Sch}
   W.Scharlau, \ {\it Quadratic and Hermitian Forms},
    Springer, Berlin, Heidelberg, New-York, Tokyo, 1985.
%
%
\bibitem{Vo4}
  V.Voevodsky, \ {\it The Milnor conjecture},
    Preprint MPI/97-8 (see also www.math.uiuc.edu/K-theory/0170/).
\end{thebibliography}
\end{document}